\documentclass[twoside,notitlepage,12pt]{article}

\usepackage{amssymb}
\usepackage[leqno]{amsmath}
\usepackage{amsfonts}
\usepackage{amsopn}
\usepackage{amstext}
\usepackage{amsthm}

\usepackage[colorlinks]{hyperref}

\renewcommand{\Bbb}{\mathbb}

\newenvironment{pf}{\begin{proof}}{\end{proof}}

\newcommand{\eps}{\varepsilon}
\renewcommand{\phi}{\varphi}
\renewcommand{\rho}{\varrho}

\newcommand{\G}{{\mathbb G}}

\newcommand{\rest}{\restriction}

\newcommand{\Ntr}{n\in{\Bbb{N}}}
\newcommand{\loe}{\leq}
\newcommand{\goe}{\geqslant}

\newcommand{\subs}{\subseteq}
\newcommand{\sups}{\supseteq}

\newcommand{\ovr}{\overline}

\newcommand{\by}{/}

\newtheorem{tw}{Theorem}[section]
\newtheorem{wn}[tw]{Corollary}
\newtheorem{lm}[tw]{Lemma}

\theoremstyle{definition}

\theoremstyle{remark}

\newcommand{\setof}[2]{\{#1\colon #2\}}

\newcommand{\sett}[2]{\{#1\}_{#2}}
\newcommand{\sn}[1]{\{#1\}} 
\newcommand{\map}[3]{#1\colon #2 \to #3} 
\newcommand{\img}[2]{#1[#2]} 

\newcommand{\ciag}[1]{{\sett{{#1}_n}{\Ntr}}}

\newcommand{\anorm}{\|\cdot\|}
\newcommand{\norm}[1]{\|#1\|}

\newcommand{\cmp}{\circ} 

\title{A proof of uniqueness of the Gurari{\u\i} space }
\author{
Wies{\l}aw Kubi\'s \footnote{Research of Kubi\'s supported in part
by the Grant IAA 100 190 901 and by the Institutional Research
Plan of the Academy of Sciences of Czech Republic No. AVOZ 101 905
03.}\; and\; S\l awomir Solecki \footnote{Research of Solecki
supported by NSF grant DMS-1001623.} }

\date{}

\begin{document}

\maketitle

\begin{abstract}
We present a short and elementary proof of isometric uniqueness of
the Gurari{\u\i} space.
\end{abstract}

\section{Introduction}

A \emph{Gurari{\u\i} space}, constructed by
Gurari{\u\i}~\cite{Gurarii} in 1965, is a separable Banach space
$\G$ satisfying the following condition: given finite dimensional
Banach spaces $X\subs Y$, given $\eps > 0$, and given an isometric
linear embedding $\map fX\G$ there exists an injective linear
operator $\map gY\G$ extending $f$ and satisfying $\norm g \cdot
\norm {g^{-1}} < 1 + \eps$. It is not hard to prove straight from
this definition that such a space is unique up to isomorphism of
norm arbitrarily close to one. The question whether the
Gurari{\u\i} space is unique up to isometry remained open for some
time. It was answered affirmatively by Lusky~\cite{Lusky} in 1976
using deep techniques developed by Lazar and
Lindenstrauss~\cite{LazLin}. Subsequently, another proof of
uniqueness was given by Henson using model theoretic methods of
continuous logic. (This proof remains unpublished.) The natural
question whether there is an elementary proof of uniqueness
occurred to several mathematicians. This question was made current
by recent increased interest in universal homogeneous structures
and their automorphism groups; see, for example, \cite{KPT} and
\cite{Mel}. The aim of this note is to provide just such a simple
and elementary proof of isometric uniqueness of the Gurari{\u\i}
space. This proof is given in Section~\ref{S:uni}. In
Section~\ref{S:univ}, we give an elementary proof showing
isometric universality of the Gurari{\u\i} space among separable
Banach spaces. Our argument uses only the basic Gurari{\u\i}
property.

In order to state the theorem precisely, we introduce some
notions. Let $X$, $Y$ be Banach spaces, $\eps>0$. A linear
operator $\map fXY$ is an \emph{$\eps$-isometry} if for $x\in X$
with $\norm{x}=1$
$$ (1+\eps)^{-1} < \norm {f(x)}|< 1+\eps.$$ We use strict inequalities for
the sake of convenience. In particular, in the case of finite
dimensional spaces, every $\eps$-isometry is an $\eps'$-isometry
for some $0 < \eps' < \eps$. Note that the inverse of a bijective
$\eps$-isometry is again an $\eps$-isometry. By an \emph{isometry}
we mean a linear operator $\map fXY$ that is an $\eps$-isometry
for every $\eps>0$, that is, $\norm {f(x)} = \norm x$ holds for
every $x\in X$. (A word of caution about our terminology may be in
place: in the literature, such functions are often called
\emph{isometric embeddings}, with the word ``isometry" reserved
for a \emph{bijective} isometric embedding.)

We will give a proof of the following theorem.

\begin{tw}\label{tmejn}
Let $E$, $F$ be separable Gurari{\u\i} spaces, $\eps>0$. Assume
$X\subs E$ is a finite dimensional space and $\map fXF$ is an
$\eps$-isometry. Then there exists a bijective isometry $\map hEF$
such that $\norm {h\rest X - f} < \eps$.
\end{tw}

By taking $X$ to be the trivial space, we obtain the following
corollary.

\begin{wn}[Lusky~\cite{Lusky}]
The Gurari{\u\i} space is unique up to a bijective isometry.
\end{wn}

\section{Proof of uniqueness of the Gurari{\u\i} space}\label{S:uni}

\begin{lm}\label{llawofreturn}
Let $X$, $Y$ be finite dimensional Banach spaces and let $\map
fXY$ be an $\eps$-isometry, for some $\eps>0$. There exist a finite
dimensional Banach space $Z$ and isometries $\map iXZ$ and $\map
jYZ$ such that $\norm { j \cmp f - i } \leq \eps$.
\end{lm}

\begin{pf}
By $\anorm_X$, $\anorm_Y$ we denote the norms of $X$ and $Y$,
respectively.

First assume that $f$ is onto, that is, $f[X]=Y$. Consider the
space $X\oplus Y$ and the canonical embeddings $\map iX{X \oplus
Y}$ and $\map j{Y}{ X \oplus Y}$. We aim to define a suitable norm
$\anorm'$ on $X\oplus Y$. With $x^*\in S_X^*$ associate the
functional ${\ovr x^*} = x^* f^{-1}$ on $Y$. Define
$$\phi_X(x,y) = \sup_{x^* \in S_X^*} \left| x^*(x) + \frac1{\norm{\ovr x^*}_Y^*} \ovr x^*(y)\right|.$$
It is clear that $\phi_X$ is a seminorm on $X \oplus Y$. Observe
that $\phi_X(x,0) = \norm{x}_X$ and $\phi_X(0,y) \loe \norm{y}_Y$.
Interchanging the roles of $X$ and $Y$ and of $f\colon X\to Y$ and
$f^{-1}\colon Y\to X$ define in an analogous way $\phi_Y$ so that
it is a seminorm on $X \oplus Y$ such that $\phi_Y(x,0) \loe
\norm{x}_X$ and $\phi_Y(0,y) = \norm{y}_Y$. Finally, define
$$\norm {(x,y)}' = \max \Bigl\{ \phi_X(x,y), \phi_Y(x,y),  \eps_1\norm{x}_X,
\eps_1 \norm{y}_Y \Bigr\},$$ where
$\eps_1 = \frac{\eps}{1+\eps}$.
Now $\anorm'$ is a norm on $X \oplus
Y$ and, since $\eps_1 \leq 1$, we have that $\norm{(x,0)}' =
\norm{x}_X$ and $\norm{(0,y)}' = \norm{y}_Y$. Hence, $i$ and $j$
are isometries.

We check that $\norm{jf(x) - i(x)}' \leq \eps$ for $x\in S_X$. Set
$u = jf(x) - i(x) = (-x, f(x))$. From the inequalities
$(1+\eps)^{-1} < \norm{ x^*f^{-1} }_{Y}^* < 1+\eps$ for $x^*\in
S_X^*$, we obtain
\begin{align*}
\phi_X(u) &= \sup_{x^*\in S_X^*}\left| x^*(-x) + \frac1{\norm{\ovr
x^*}_Y^*} \ovr x^*(f(x)) \right|
= \sup_{x^*\in S_X^*} \left(\left|  \frac1{\norm{x^* f^{-1}}_{Y}^*} - 1 \right|\cdot |x^*(x)| \right) \\
&\leq \sup_{x^*\in S_X^*} \left|  \frac1{\norm{x^* f^{-1}}_{Y}^*}
- 1 \right| \leq \eps.
\end{align*}
Similarly, from $(1+\eps)^{-1} < \norm{ y^*f }_{X}^* < 1+\eps$ for
$y^*\in S_Y^*$, we get
\begin{align*}
\phi_Y(u) &= \sup_{y^*\in S_Y^*}\left| y^*(f(x)) +
\frac1{\norm{\ovr y^*}_X^*} \ovr y^*(-x) \right|
= \sup_{y^*\in S_Y^*} \left(\left| 1- \frac1{\norm{y^* f}_{X}^*}  \right|\cdot |y^*(f(x))| \right) \\
&\leq \sup_{y^*\in S_Y^*} \left(\left| 1- \frac1{\norm{y^*
f}_{X}^*} \right|\cdot \norm{y^*f}_X^* \right)= \sup_{y^*\in
S_Y^*} \left| \norm{y^* f}_{X}^*-1 \right| \leq \eps.
\end{align*}
Finally, since
\[
\eps_1 \norm{-x}_X \leq \eps\; \hbox{ and }\; \eps_1
\norm{f(x)}_Y\leq \frac{\eps}{1+\eps}(1+\eps) =\eps,
\]
we conclude that ${\norm u}' \leq \eps$, as required.

Now we consider the general case when $f$ is not necessarily onto.
The conclusion above gives a norm $\anorm'$ on $X\oplus f[X]$.
Take $(X\oplus f[X])\oplus Y$ regarded as the $\ell_1$ sum of the
Banach spaces $(X\oplus f[X],\, \anorm')$ and $(Y,\, \anorm_Y)$.
Pass to the quotient Banach space
\[
Z = ((X\oplus f[X])\oplus Y)/\{ (0, f(v), -f(v))\colon v\in X\}
\]
with the quotient norm and with the canonical embeddings of $X$
and $Y$. This space is as required. Note that $Z$ is canonically
isometric to $X\oplus Y$ equipped with the norm
\[
\norm{(x,y)} = \inf_{v\in X} \Bigl(\norm{(x,f(v))}' +
\norm{y-f(v))}_Y\Bigr). \qedhere
\]
\end{pf}

\begin{lm}\label{lgherik}
Let $E$ be a Gurari{\u\i} space and let $\map fXY$ be an
$\eps$-isometry, where $X$ is a finite dimensional subspace of $E$
and $\eps>0$. Then for every $\delta > 0$ there exists a
$\delta$-isometry $\map gYE$ such that $\norm { g\cmp f - {\rm
id}_X } < \eps$.
\end{lm}

\begin{pf}
Choose $0 < \eps' < \eps$ so that $f$ is an $\eps'$-isometry.
Choose $0 < \delta' < \delta$ such that $(1 + \delta')\eps' <
\eps$. By Lemma~\ref{llawofreturn}, there exist a finite
dimensional space $Z$ and isometries $\map iXZ$ and $\map jYZ$
satisfying $\norm{j\cmp f - i} \leq \eps'$. Since $E$ is
Gurari{\u\i}, there exists a $\delta'$-isometry $\map hZE$ such
that $hj(x)=x$ for $x\in X$. Let $g = h\cmp j$. Clearly, $g$ is a
$\delta$-isometry. Finally, given $x\in S_X$, we have
$$\norm{ gf(x) - x } = \norm {hjf(x) - hi(x)} < (1+\delta') \norm{ jf(x) - i(x)}
\leq (1+\delta')\eps' < \eps,$$ as required.
\end{pf}

\begin{pf}[Proof of Theorem~\ref{tmejn}]
Fix a sequence $\ciag \eps$ of positive real numbers. The precise
conditions on $\ciag \eps$ will be specified later. We define
inductively sequences of linear operators $\ciag{f}$, $\ciag g$
and finite dimensional subspaces $\ciag{X}$, $\ciag{Y}$ of $E$ and
$F$, respectively, so that the following conditions are satisfied:
\begin{enumerate}
    \item[(0)] $X_0 = X$, $Y_0 = \img fX$, and $f_0 = f$;
    \item[(1)] $\map {f_n} {X_{n}} {Y_{n}}$ is an $\eps_{n}$-isometry;
    \item[(2)] $\map {g_n} {Y_{n}} {X_{n+1}}$ is an $\eps_{n+1}$-isometry;
    \item[(3)] $\norm{g_n f_n(x) - x} \loe \eps_{n} \norm x$ for $x\in X_{n}$;
    \item[(4)] $\norm{f_{n+1}g_n(y) - y} \loe \eps_{n+1} \norm y$ for $y\in
    Y_{n}$;
    \item[(5)] $X_n\subseteq X_{n+1}$, $Y_n\subseteq Y_{n+1}$,
    $\bigcup_nX_n$ and $\bigcup_nY_n$ are dense in $E$ and $F$,
    respectively.
\end{enumerate}
Condition (0) tells us how to start the inductive construction. We
pick $\eps_0>0$ so that (1) holds for $n=0$ and $\eps_0< \eps$.
Suppose $f_i$, $X_i$, $Y_i$, for $i\leq n$, and $g_i$, for $i <
n$, have been constructed. We easily find $g_n$, $X_{n+1}$,
$f_{n+1}$ and $Y_{n+1}$, in this order, using Lemma~\ref{lgherik}.
Condition (5) can be secured by defining $X_{n+1}$ and $Y_{n+1}$
to be appropriately enlarged $g_n[Y_n]$ and $f_{n+1}[X_{n+1}]$,
respectively. Thus, the construction can be carried out.

Fix $\Ntr$ and $x\in X_{n}$ with $\norm x = 1$. Using (4) and (1),
we get
$$\norm{ f_{n+1} g_n f_n(x) - f_n(x) } \loe \eps_{n+1} \norm{f_n(x)} \loe \eps_{n+1}(1 + \eps_{n}).$$
Using (1) and (3), we get
$$\norm{ f_{n+1} g_n f_n(x) - f_{n+1}(x) } \loe \norm{f_{n+1}} \cdot \norm{ g_n f_n(x) - x} \loe (1 + \eps_{n+1})\cdot
\eps_{n}.$$ These inequalities give
\begin{equation}
\norm{ f_n(x) - f_{n+1}(x) } \loe \eps_{n} + 2\eps_n \eps_{n+1} +
\eps_{n+1}. \tag{$\dagger$}\label{eqsztylet}
\end{equation}
Now it is clear that if the series $\sum_{\Ntr}\eps_n$ converges,
then the sequence $\sett{f_n(x)}{\Ntr}$ is Cauchy. Let us make a
stronger assumption, namely that
\begin{equation}
2\eps_0\eps_1+\eps_1+\sum_{n=1}^\infty (\eps_{n} +
2\eps_n\eps_{n+1} + \eps_{n+1}) < \eps - \eps_0.
\tag{$\ddagger$}\label{eqssztilett}
\end{equation}
Given $x\in \bigcup_{\Ntr}X_n$, define $h(x) = \lim_{n \goe
m}f_n(x)$, where $m$ is such that $x\in X_{m}$. Then $h$ is an
$\eps_n$-isometry for every $\Ntr$, hence it is an isometry.
Consequently, it uniquely extends to an isometry on $E$, which we
denote also by $h$. Furthermore, (\ref{eqsztylet}) and
(\ref{eqssztilett}) give
$$\norm{f(x) - h(x)} \loe
\sum_{n=0}^\infty \eps_{n} + 2\eps_n\eps_{n+1} + \eps_{n+1} <
\eps.$$

It remains to see that $h$ is a bijection. To this end, we check
as before that $\sett{g_n(y)}{n\goe m}$ is a Cauchy sequence for
every $y\in Y_{m}$. Once this is done, we obtain an isometry
$g_\infty$ defined on $F$. Conditions (3) and (4) tell us that
$g_\infty \cmp h = {\rm id}_E$ and $h \cmp g_\infty = {\rm id}_F$,
and the proof is complete.
\end{pf}

\section{On universality of the Gurari{\u\i} space}\label{S:univ}

It is known that the Gurari{\u\i} space is isometrically universal
among separable Banach spaces. Indeed, as pointed out by
Gevorkjan~\cite{Gevork}, universality follows from the results of
Lazar and Lindenstrauss~\cite{LazLin} and Michael and
Pe{\l}czy\'nski~\cite{MP}: the dual of the Gurari{\u\i} space is a
non-separable $L_1$ space, therefore the Gurari{\u\i} space
contains an isometric copy of $C([0,1])$. The reader may also
consult the recent paper \cite{Avi} for another approach.

We conclude with applying our method to proving universality
directly, without referring to the structure of the dual or to
universality of other Banach spaces.

\begin{lm}\label{ltgereq}
Let $X_0, X_1, Y_0$ be finite dimensional Banach spaces such that
$X_0 \subs X_1$ and let $\map f{X_0}{Y_0}$ be an $\eps$-isometry,
where $\eps > 0$. Then there exist a finite dimensional Banach
space $Y_1$ containing $Y_0$ and an isometry $\map g{X_1}{Y_1}$
such that
$$\norm { g\rest X_0 - f } < \eps.$$
\end{lm}

\begin{pf}
A standard and well known amalgamation property for Banach spaces (already used in the proof of Lemma~\ref{llawofreturn} above)
says that there exist $W\sups Y_0$ and an $\eps$-isometry $\map
{f'}{X_1}W$ such that $f' \rest X_0 = f$. More precisely, $W =
(X_1 \oplus Y_0)\by \Delta$, where $X_1\oplus Y_0$ is endowed with
the $\ell_1$-norm and
$$\Delta = \setof{(z,-f(z))}{z\in X_0}.$$
The space $Y_0$ is naturally identified with the subspace of $W$
and $f'(x)$ is the equivalence class of $(x,0)$ (where $x\in
X_1$).

Finally, the desired isometry $g$ is provided by
Lemma~\ref{llawofreturn}.
\end{pf}

\begin{tw}
Every separable Banach space can be isometrically embedded into
the Gurari{\u\i} space.
\end{tw}

\begin{pf}
Let $\G$ denote the Gurari{\u\i} space. Fix a separable Banach
space $X$ and let $\ciag X$ be a chain of finite dimensional
spaces such that $X_0=\sn0$ and $\bigcup_{\Ntr}X_n$ is dense in
$X$. In case $X$ is finite dimensional, we set $X_n = X$ for $n >
0$. We inductively define $\map {f_n}{X_n}\G$ so that
\begin{enumerate}
    \item[(i)] $f_n$ is a $2^{-n}$-isometry,
    \item[(ii)] $\norm{ f_{n+1} \rest X_n - f_n } < 2\cdot 2^{-n}$,
\end{enumerate}
for every $\Ntr$. We set $f_0 = 0$. Suppose $f_n$ has already been
defined. Let $Y = \img {f_n}{X_n}$. Using Lemma~\ref{ltgereq}, we
find a finite dimensional space $W \sups Y$ and an isometry $\map
g{X_{n+1}}W$ such that $\norm {g\rest X_n - f_n} < 2^{-n}$. Using
the property of the Gurari{\u\i} space, we find a
$2^{-(n+1)}$-isometry $\map hW\G$ such that $h \rest Y$ is the
inclusion $Y\subs \G$. Now set $f_{n+1} = h\cmp g$. Given $x\in
X_n$ with $\norm x = 1$, we have that $\norm{g(x) - f_n(x)} <
2^{-n}$ and hence
$$\norm{f_{n+1}(x) - f_n(x)} = \norm{h(g(x)) - h(f_n(x))} < (1 + 2^{-(n+1)})\cdot 2^{-n} \loe 2\cdot 2^{-n}.$$
This shows (ii). Finally, we obtain a sequence $\ciag f$ that is
pointwise Cauchy on each $X_n$. By (i) and (ii), $f_\infty(x) =
\lim_{n\to\infty}f_n(x)$ is a well-defined linear isometry on
$\bigcup_{\Ntr}X_n$. This isometry extends uniquely to an isometry
$\map fX\G$.
\end{pf}

\medskip

{\bf Acknowledgement.} We thank Ward Henson and Julien Melleray
for their comments on an earlier version of this paper.

{\small Kubi{\'s}'s address:}\\
{\small Mathematical Institute, Academy of Sciences of the Czech Republic, Prague, Czech Republic}\\
{\small Institute of Mathematics,}
{\small Jan Kochanowski University, Kielce, Poland}\\
{\small\texttt{kubis@math.cas.cz}, \texttt{wkubis@pu.kielce.pl}

\medskip
{\small Solecki's address:}\\
{\small Department of Mathematics, University of Illinois, Urbana, Illinois 61801, USA}\\
{\small Institute of Mathematics, Polish Academy of Sciences, Warsaw, Poland}\\
{\small\texttt{ssolecki@math.uiuc.edu}}

\end{document}